\documentclass{mcom-l}
\usepackage{psfig}
\newtheorem{theorem}{Theorem}[section]

\newtheorem{algorithm}[theorem]{Algorithm}
\theoremstyle{definition}

\theoremstyle{remark}
\newtheorem{remark}[theorem]{Remark}

\numberwithin{equation}{section}
\newcommand{\be}{\begin{equation}}
\newcommand{\ee}{\end{equation}}
\newcommand{\ea}{\end{array}}

\newcommand{\Ai}{\mbox{Ai}}
\newcommand{\Bi}{\mbox{Bi}}
\newcommand{\Gi}{\mbox{Gi}}
\newcommand{\Hi}{\mbox{Hi}}
\newcommand{\phase}{\mbox{ph}}

\begin{document}

\title[Non-Oscillating Integrals for the 
Scorer functions]
{On Non-Oscillating Integrals for Computing 
Inhomogeneous Airy Functions}

\author{Amparo Gil, Javier Segura}
\address{
  Instituto de Bioingenier\'{\i}a,
  Universidad Miguel Hern\'andez, Edificio La
  Galia. 03202-Elche (Alicante), Spain}
\email{amparo@titan.ific.uv.es}
\email{segura@flamenco.ific.uv.es}

\author{Nico M. Temme}
\address{CWI,  
P.O. Box 94079, 
1090 GB Amsterdam, 
The Netherlands}
\email{nicot@cwi.nl}


\begin{abstract}

Integral representations are considered of solutions of the inhomogeneous
Airy differential equation $w''-z\,w=\pm1/\pi$. 
The solutions of these equations
are also known as Scorer functions. Certain functional
relations for these functions are used to confine the discussion to
one function and to a certain sector in the complex plane. By using steepest
descent methods from  asymptotics, the standard integral representations of the
Scorer functions are modified in order to obtain non-oscillating integrals for
complex values of $z$. In this way stable representations for numerical
evaluations of the functions are obtained. 
The methods are illustrated
with numerical results.

\end{abstract}

\maketitle

\section{Introduction}

Airy functions are solutions of the differential equation

\begin{equation}
\frac{d^2\,w}{d\,z^2}-z\,w=0.\label{i1}
\end{equation}

Two linearly independent solutions that are real for real values of $z$ are
denoted by $\Ai (z)$ and $\Bi (z)$. They have the integral representations

\begin{equation}
\begin{array}{ll}
\Ai (z)&=\frac1\pi\,\int_{0}^{\infty}\,\cos\left(zt+\frac13t^3\right)\,dt ,\\
\\
\Bi (z)&=\frac1\pi\,\int_{0}^{\infty}\,\sin\left(zt+\frac13t^3\right)\,dt+
\frac1\pi\,\int_{0}^{\infty}
e^{zt-\frac13t^3}\,dt \label{i2}
\end{array}
\end{equation}
where we assume that  $z$ is real.

In this paper we concentrate on so-called Scorer functions (\cite{Lee, Sco}), 
which are particular solutions of the
non-homogeneous Airy differential equation.
We have 

\begin{equation}
w''-z\,w=-1/\pi,\quad {\hbox{\rm with solution}}\quad
\Gi (z)=\frac1\pi\, \int_{0}^{\infty}\,\sin\left(zt+\frac13t^3\right)
\,dt,\quad 
z\in\mathbb{R},\label{i3}
\end{equation}
and 
\begin{equation}
w''-z\,w=1/\pi,\quad \hbox{\rm with solution}\quad
\Hi (z)=\frac1\pi\, \int_{0}^{\infty} e^{zt-\frac13t^3}\,dt,\quad z\in 
\mathbb{C}.
\label{i4}
\end{equation}
Initial values are

\begin{equation}
\begin{array}{ll}
\Gi (0) &=\frac{1}{2}\Hi(0) =\ \frac{1}{3}\Bi(0)=\ \ \,\frac{1}{\sqrt{{3}}}
\Ai(0)=\ \
\frac{1}{3^{7/6}\Gamma(\frac{2}{3})},\\
 \\
\Gi '(0)&=\frac{1}{2}\Hi'(0)=\frac{1}{3}\Bi'(0) =-\frac{1}{\sqrt{{3}}}\Ai'(0) 
=\frac{1}{3^{5/6}\Gamma(\frac{1}{3}). }\label{i5}
\end{array}
\end{equation}

From (\ref{i2}), (\ref{i3}) and (\ref{i4}) it follows that
\begin{equation}
\Gi (z)+\Hi (z)=\Bi (z).\label{i6}
\end{equation} 

In the next section we give contour integrals from which representations of
$\Ai (z),\mbox{Bi}(z)$ and 
$\Gi (z)$ follow for complex values of $z$. 
Just like  $\mbox{Ai}(z)$ and $\mbox{Bi}(z)$, the Scorer functions 
$\Gi (z)$ and
$\Hi (z)$ are entire functions. 

A survey on computational aspects of special functions, including
information on Airy functions, can be found in \cite{Loz}; \cite{Nis}
has a public web site that includes an extensive treatment of
Scorer functions. For complex
values of $z$ the Airy functions are available in the  Bessel function
algorithms of \cite{Amo}; see also \cite{Cor} and \cite{Sch}. Computer 
algebra
systems as Maple and Mathematica also have Airy functions available.
The Scorer functions are considered in \cite{Leo}, where coefficients of
Chebyshev  expansions are given for real $z$. Asymptotic expansions for
$\Hi(z)$ are given  in \cite{Ext} and \cite{Olv}. The paper by Scorer
(\cite{Sco}) gives tables to 7 decimals of $\Hi(-z)$ and $\Gi(z)$ for $0\le
z\le 10$.

Efficient algorithms for computing the Scorer functions
 in restricted domains of
the complex plane can be based on Maclaurin series and  asymptotic series.
These  domains can be bridged by using the differential equations or the
integral representations. 

In \cite{Fab} methods were presented based on the
differential equations,  which are set up as boundary-value methods, providing
stable algorithms for all values of $z$. 

The purpose of this  paper is to give stable integral representations
for $\Gi(z)$ and $\Hi(z)$.  We modify the integrals in (\ref{i3}) and
(\ref{i4}) by using methods from asymptotics. As in \cite{Fab}, the resulting
integrals can be used for any value of the complex parameter $z$. 
We also indicate how similar methods can be used for the Airy function itself.

\section{Asymptotic properties of the Airy and Scorer functions}
We need a few properties of the Airy and Scorer functions. More information can
be found in \cite{Abr}, \cite{Nis},  \cite{Olv} and \cite{Tem2}. In particular 
\cite{Olv}, 
Chapter 11, discusses numerically satisfactory solutions of the differential 
equations for $\Gi(z)$ and $\Hi(z)$. The asymptotic properties 
of the Airy and Scorer functions are important in this discussion.

We write, as in \cite{Olv}, 
\begin{equation}
\Ai_0(z)=\Ai(z),\quad \Ai_1(z)=\Ai\left(e^{-2\pi i/3}z\right),
\quad \Ai_{-1}(z)=\Ai\left(e^{2\pi i/3}z\right). \label{t1}
\end{equation}
We have the representations 
\begin{equation}
\Ai_j(z)=\frac{e^{2j\pi i/3}}{2\pi i}
\int_{\mathcal{C}_{-j}}
e^{-zt+\frac13t^3}\,dt,\quad j=0,\pm1,\label{t2}
\end{equation}
where the contours ${\mathcal{C}}_j$ are given in Figure 2.1. 
Because 
$\int_{{\mathcal{C}}_0\cup{\mathcal{C}}_{1}\cup{\mathcal{C}}_{-1}}e^{-
zt+\frac13t^3}\,dt=0$,
we have the following linear combination of three solutions of (\ref{i1}):
\begin{equation}
\Ai(z)+e^{-2\pi i/3} \Ai_1(z)+e^{2\pi i/3} \Ai_{-1}(z)=0.\label{t3}
\end{equation}

\vspace*{0.4cm}
\centerline{\protect\hbox{\psfig{file=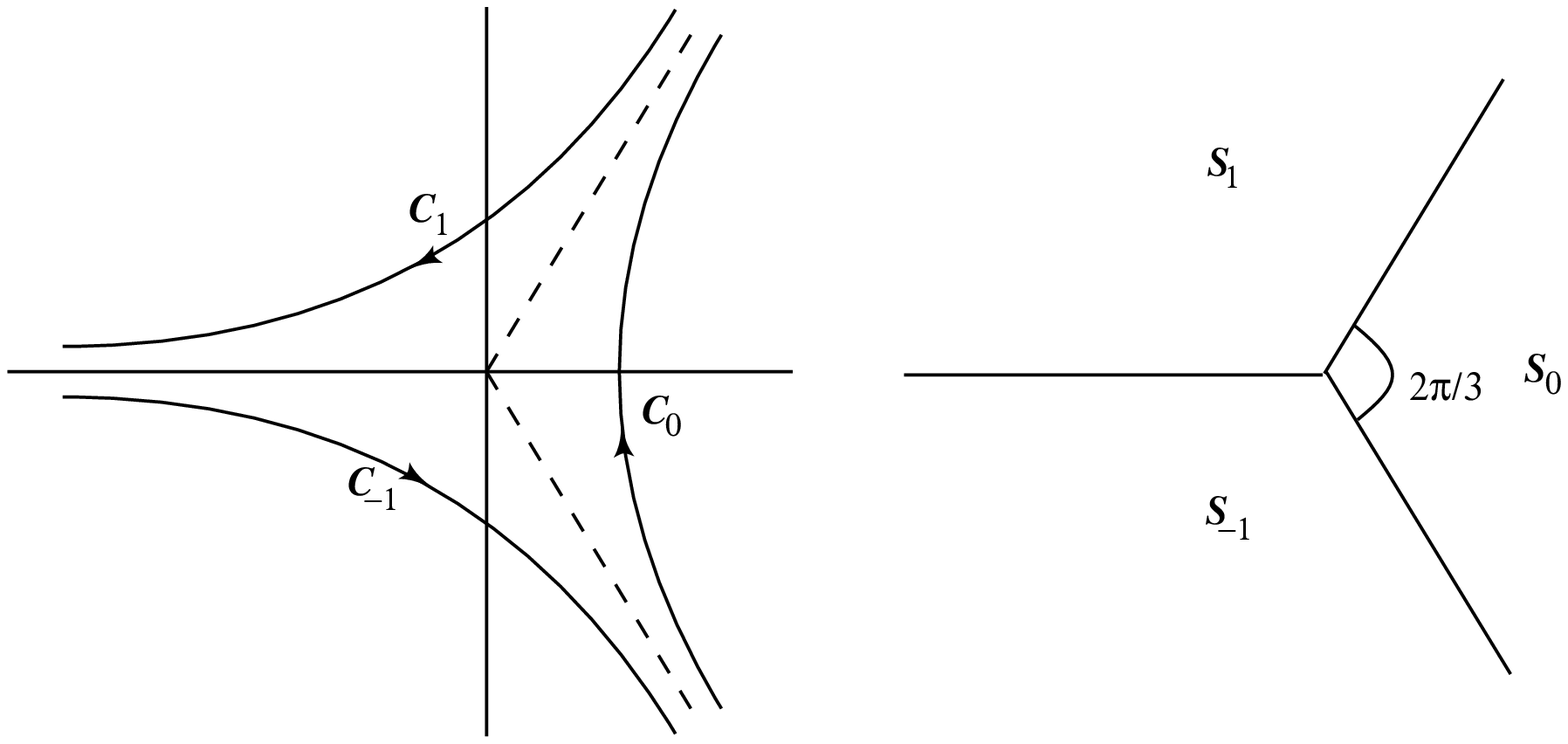,width=10.cm}}}
\noindent	
{\bf Figure 2.1.}\quad Three contours ${\mathcal{C}}_j$ of integration
 for the Airy
integrals  in (\ref{t2}), and sectors $S_j$ where $\Ai_j(z)$ are recessive.

\vspace*{0.4cm}

The first integral in (\ref{i2}) follows from deforming
 the contour ${{{\mathcal 
{C}}}_0}$ 
in (\ref{t2}) into the imaginary axis.  
The function $\Bi(z)$ can be written as
\begin{equation}
\Bi(z)=e^{\pi i/6} \Ai_{-1}(z)+e^{-\pi i/6} \Ai_1(z)\label{t4}
\end{equation}
and  the second representation in (\ref{i2}) 
follows by deforming contour ${{\mathcal{C}}_1}$ into the
 positive imaginary axis
and $(-\infty,0]$,
and the contour ${{\mathcal{C}}_{-1}}$  into
$(-\infty,0]$ and the negative imaginary axis.

The function $\Ai_j(z)$ is recessive at infinity in the sector
 $S_j, j=0, \pm1$, 
the function being exponentially small at infinity along any ray interior to
this sector. On the other hand,  $\Ai_j(z)$ is dominant at 
infinity in the sectors
$S_{j-1}$ and $S_{j+1}$ (the suffix $j$ is enumerated modulo 3), and is
exponentially large at infinity along the rays interior to  these sectors. 
$\Bi(z)$ is dominant at infinity within all three sectors $S_j$. A pair of Airy
functions comprises a numerically satisfactory pair at infinity  within a
sector if only one function is dominant. For example, the pair 
$\{\Ai(z), \Bi(z)\}$ comprises such a pair only in $S_0$ (and on
the negative real axis, where none of the two is dominant, but where the phases
in their oscillations differ by $\frac12\pi$).

If one of the functions  $\Gi(z), \Hi(z)$ is computed, we may use (\ref{i6}) 
to compute the other one (we assume in this paper that $\Bi(z)$ and all 
other Airy  functions are available), but we need to know if (1.6) is
numerically stable for that computation. For example, because
$\Gi(z)$ is only of algebraic growth in $S_0$ (as we will see soon), we cannot
compute $\Gi(z)$ from (\ref{i6}) in $S_0$.

For the Scorer functions we have the following asymptotic expansions (cf.
\cite{Olv}, 431--432):
\begin{equation}
\Hi(z)\sim-\frac1{\pi 
z}\left[1+\frac1{z^3}\sum_{s=0}^\infty\frac{(3s+2)!}{s!(3z^3)^s}\right],
\quad z\to\infty,\quad |\phase(-z)|\le\frac23\pi-\delta,\label{t5}
\end{equation}
\begin{equation}
\Gi(z)\sim \frac1{\pi 
z}\left[1+\frac1{z^3}\sum_{s=0}^\infty\frac{(3s+2)!}{s!(3z^3)^s}\right],
\quad z\to\infty,\quad |\phase\ z|\le\frac13\pi-\delta,\label{t6}
\end{equation}
$\delta$ being an arbitrary positive constant. For (\ref{t6}) the domain 
for $\phase\ z$ is
not given in \cite{Olv}, but it follows from the same
method mentioned for $\Hi(z)$ in \cite{Olv}, p. 432. In other parts of the
complex plane we cannot represent the Scorer functions
 by a single expansion with
leading term ${\mathcal{O}}(1/z)$. 

From the results in (\ref{t5}) and (\ref{t6}) and the dominant asymptotic 
behavior
of $\Bi(z)$ in all sectors $S_j, j=0,\pm1$, we conclude that (cf. (\ref{i6}))
$\Bi(z)$ is a dominant term for $\Hi(z)$ in $S_0$ and for $\Gi(z)$ in $S_1\cup
S_{-1}$. It follows that we need algorithms for the computation of $\Gi(z)$ for
$z\in S_0$ and for $\Hi(z)$ for $z\in S_1\cup S_{-1}$ (where the asymptotic
expansions  (\ref{t5}) and (\ref{t6}) are valid).
 The relation in (\ref{i6}) can be used
for computing the functions in the complements of these domains (where the
functions have the dominant behaviour of $\Bi(z)$).

A further reduction of domains follows from the connection formula
\footnote{With thanks to the referee.}
\begin{equation}
\Hi(z)=e^{\pm2\pi i/3}\Hi\left(ze^{\pm2\pi i/3}\right)+2e^{\mp\pi 
i/6}\Ai\left(ze^{\mp2\pi
i/3}\right).\label{t7}
\end{equation}
To prove this relation observe that the first term in the
right-hand side satisfies the differential equation for $\Hi(z)$, and that,
hence, that term can be written as a linear combination of $\Hi(z)$ and
solutions of the homogeneous equation; the initial values
 in (\ref{i5}) can be used to identify these solutions. 

For example, we can use (\ref{t7}) with the upper signs for $z$ in the sector
$\frac13\pi<\phase\ z<\frac23\pi$. Then, $\Hi(z)$ can be expressed in terms of 
$\Hi$
in the sector
$-\pi<\phase\ z<-\frac23\pi$ plus an Airy function in the sector
$-\pi/3<\phase\ z<0$. 
We see that both functions in the right-hand side
 of (\ref{t7}) are not dominant 
in the respective sectors, and, hence, this representation is stable. 

A similar connection formula for $\Gi(z)$ reads
$$\Gi(z)=e^{\pm2\pi i/3}\Gi\left(ze^{\pm2\pi i/3}\right)+e^{\mp\pi 
i/6}\Ai\left(ze^{\mp2\pi
i/3}\right).$$
This formula is of no use in the sector $S_0$ because both functions in the
right-hand side are dominant, whereas $\Gi(z)$ is of algebraic growth at
infinity within $S_0$. A better formula for $z\in S_0$ follows from combining
(\ref{i6}), (\ref{t4}) and (\ref{t7}) (the latter twice,
 with upper and lower signs). This gives

\begin{equation}
\Gi(z)=-\frac12\left[e^{2\pi i/3}\Hi\left(ze^{2\pi i/3}\right)+e^{-2\pi i/3}
\Hi\left(ze^{-2\pi i/3}\right)\right].\label{t8}
\end{equation} 
For $z\in S_0$ the arguments of the $\Hi-$functions are in $S_{\pm1}$,
where these functions have expansions that follow from (\ref{t5}).

Because of 
\begin{equation}
\Hi(x-iy)=\overline{\Hi(x+iy)}\quad {\rm and}\quad
\Gi(x-iy)=\overline{\Gi(x+iy)}\label{t9}
\end{equation}
we can concentrate on non-negative
values of the imaginary part $y$ of the argument $z=x+iy$.

Conclusion.\quad\rm
The principal domain of interest for the Scorer functions is
the sector $\frac23\pi\le\phase\ z\le\pi$, where  
we concentrate on $\Hi(z)$. 
For $z$ in other sectors, and for $\Gi(z)$, 
the relations (\ref{i6}), (\ref{t7}) and (\ref{t8}) are numerically stable
for the particular cases. 

For a summary of the results of this section and algorithms we
refer to \mbox{Section 4}.

\section{The construction of non-oscillating integrals}%
We modify the integral in (\ref{i4}) such that
stable algorithms can be based on the new integral for $z$ in the sector
$\frac23\pi\le\phase\ z\le\pi$.

We write
\begin{equation}
\phi(t)=\frac13t^3-zt, \quad t=u+iv,\quad z=x+iy.\label{d1}
\end{equation}
Then the real and imaginary parts of $\phi(t)=\phi_r(u,v)+i\phi_i(u,v)$ are
given by
\begin{equation}
\begin{array}{ll}
\phi_r(u,v)&=\frac13u^3-uv^2-xu+yv,\\
\phi_i(u,v)&=u^2v-\frac13v^3-xv-yu.\label{d2}
\end{array}
\end{equation}

We are interested in a path in the complex $t-$plane on which
$\phi_i(u,v)$ is a constant, and the path should start at the origin,
as the integral in (\ref{i4}). Such a path is defined by the equation
\begin{equation}
u^2v-\frac13v^3-xv-yu=0,\label{d3}
\end{equation}
and we need real solutions of this equation. 

As summarized at the end of Section 2, we need to consider  $z-$values 
satisfying $2\pi/3\le\phase\ z\le\pi$. If $z<0$ we can integrate along 
the positive $t-$axis. For other values of $\phase\ z$ in the present
range the path of integration will be deformed
 into a curve ${\mathcal{C}}$  that
is defined by (\ref{d3}) and that
runs from the origin to $+\infty$; see Figure 3.1.

\vspace*{0.4cm}
\centerline{\protect\hbox{\psfig{file=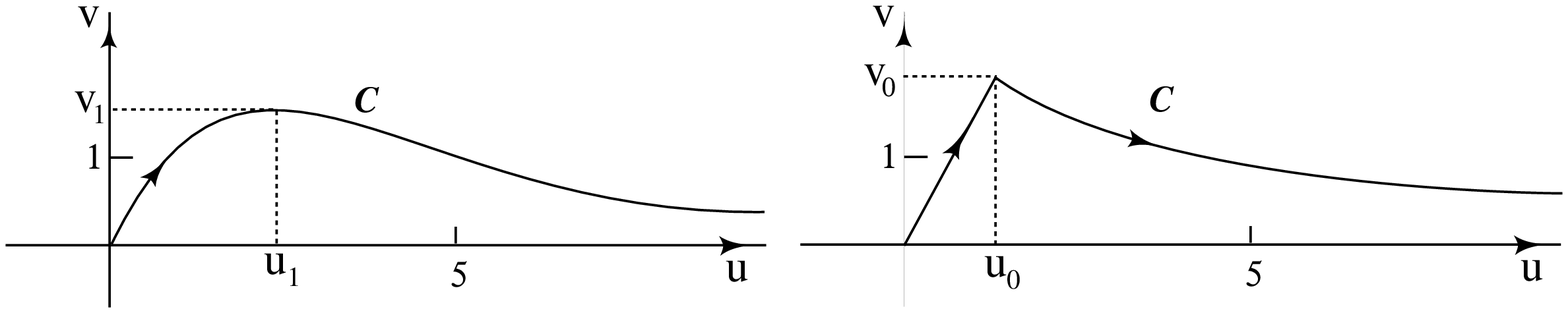,width=12.cm}}}
\noindent	
{\bf Figure 3.1.}\quad The contour ${\mathcal{C}}$ is defined by (\ref{d4}).
The left figure is for  $z$ inside the sector
 $2\pi/3\le\phase\ z\le\pi$, and the curve is defined
 by (\ref{d4}). If $\phase\ z=2\pi/3$ the curve in the right figure
 is defined by (\ref{d8}). 

\vspace*{0.4cm}

We solve the cubic equation (\ref{d3}) for $v$. The solution that passes
through the origin can be written in the form
\begin{equation}
v(u)=2\sqrt{{Q}}\sin\frac1{ 3}\theta,\label{d4}
\end{equation}
where
\begin{equation}
\theta=\arcsin\frac P{\sqrt{{Q^3}}}\in[0,\frac12\pi],\quad P=\frac32yu,\quad
Q=u^2-x,\label{d5} 
\end{equation}
with $u\ge0, x<0, 3x^2>y^2$.

To show this, we introduce $t=v/(2\sqrt{{Q}})$. Then (\ref{d3}) can be
written in the form
$$4t^3-3t=-\frac{P}{\sqrt{{Q^3}}}.$$
It is not difficult to verify that, if $u\ge0, x<0, 3x^2>y^2$, then the
modulus of the right-hand side is not larger than unity. Replacing the
left-hand side with $-\sin(3\arcsin(t))$ gives the solution in
(\ref{d5}).

It follows that
\begin{equation}
\Hi(z)=\frac1\pi\int_{0}^{\infty}\,e^{-\phi_r(u,v(u))} h(u)\,du,\label{d6}
\end{equation}
where $\phi_r(u,v)$ is given in (\ref{d2}), $v(u)$ in (\ref{d4}) and
\begin{equation}
h(u)=\frac{dt}{du}=\frac{d[u+iv(u)]}{du}=1+i\frac{dv(u)}{du}=1+
i\frac{2uv-y}{v^2-u^2+x}.\label{d7}
\end{equation}

If $\phase\ z =2\pi/3$ then $y=-x\sqrt{{3}}, and x\le0$, equation
(\ref{d3}) can be solved explicitly. The two solutions are
\begin{equation}
v=u\sqrt{{3}} \quad {\rm and} \quad u=-\frac{3x+v^2}{v\sqrt{{3}}}.\label{d8}
\end{equation}
In this case the path of integration ${\mathcal{C}}$
 runs from the origin to the 
point 
$$t_0= u_0+iv_0=\sqrt{{-x/2}}+i\sqrt{{-3x/2}}$$
along the line $v=u\sqrt{{3}}$, and for $u\ge u_0$
 the path ${\mathcal{C}}$ follows
the hyperbola defined by the second solution given in (\ref{d8}); see the right
figure in Figure 3.1. The point $t_0=u_0+iv_0=\sqrt{{z}}=\sqrt{{x-
ix\sqrt{{3}}}}, x
\le0$,  is a saddle point of the function $\phi(t)$ defined in (\ref{d1}).

We can also solve (\ref{d3}) for $u$, which gives 
\begin{equation}
u=\frac{y-R}{2v},\quad R=\sqrt{{y^2+4v^2(x+\frac13v^2)}},\label{d9}
\end{equation}
where the square root is non-negative. This solution should be used for
$0\le v\le v_1, 0\le u\le u_1$, where 
$$v_1=v_1(x,y)=\sqrt{{\frac32\left(-x-\sqrt{{x^2-y^2/3}}\right)}},$$
the smallest positive $v-$value  for which
$R=0$ and $u_1=y/(2v_1)$.  For  $0\le v\le v_1, u\ge u_1$, 
we use $u=(y+R)/(2v)$; see Figure 3.1.

When integrating with
respect to
$v$, the integral in (\ref{d6}) can be  written as
\begin{equation}
\Hi(z)=\frac1\pi\left[\int_0^{v_1}\,e^{-\phi_r(u^-(v),v)} h(v)\,dv
+\int_{v_1}^{0}\,e^{-\phi_r(u^+(v),v)} h(v)\,dv\right],\label{d10}
\end{equation}
where $\phi_r(u,v)$ is given in (\ref{d2}), $u^\pm(v)=(y\pm R)/(2v)$ (cf. 
(\ref{d9}))
and 
\begin{equation}
h(v)=\frac{dt}{dv}=\frac{d(u+iv)}{dv}=\frac{du}{dv}+i=
\frac{v^2-u^2+x}{2uv-y}+i.\label{d11}
\end{equation}

\remark
For the sector $0\le\phase\ z\le\frac23\pi $ we can use a similar method, 
although we don't need to consider this sector.  
If $y^2-3x^2\ge0$ the quantity $R$ of (\ref{d9}) is defined for all values of 
$v$, and the first equation in (\ref{d9}) defines a path
 ${\mathcal{L}}$ going from the origin 
to $\infty\exp(2\pi i/3)$. Because there is a hill at $\infty\exp(\pi i/3)$, 
we cannot replace the integration path in (\ref{i4}) by ${\mathcal{L}}$. 
We need an extra integral from $\infty\exp(2\pi i/3)$ to $+\infty$, and that 
integral
gives an Airy function; see Figure 3.2. In this way we obtain 
\begin{equation}
\Hi(z)=\frac1\pi\int_{{\mathcal{L}}}\,e^{-\phi_r(u,v)} h(v)\,dv+
2e^{-\pi i/6}\Ai\left(ze^{-2\pi i/3}\right),\label{d12}
\end{equation}
with $h(v)$ given in (\ref{d11}) and the relation between $u$ and $v$ 
given in (\ref{d9}).  We see that  the first term in the right-hand 
side of (\ref{t7}) with upper signs corresponds with the
 integral in (\ref{d12}). 

\vspace*{0.4cm}
\centerline{\protect\hbox{\psfig{file=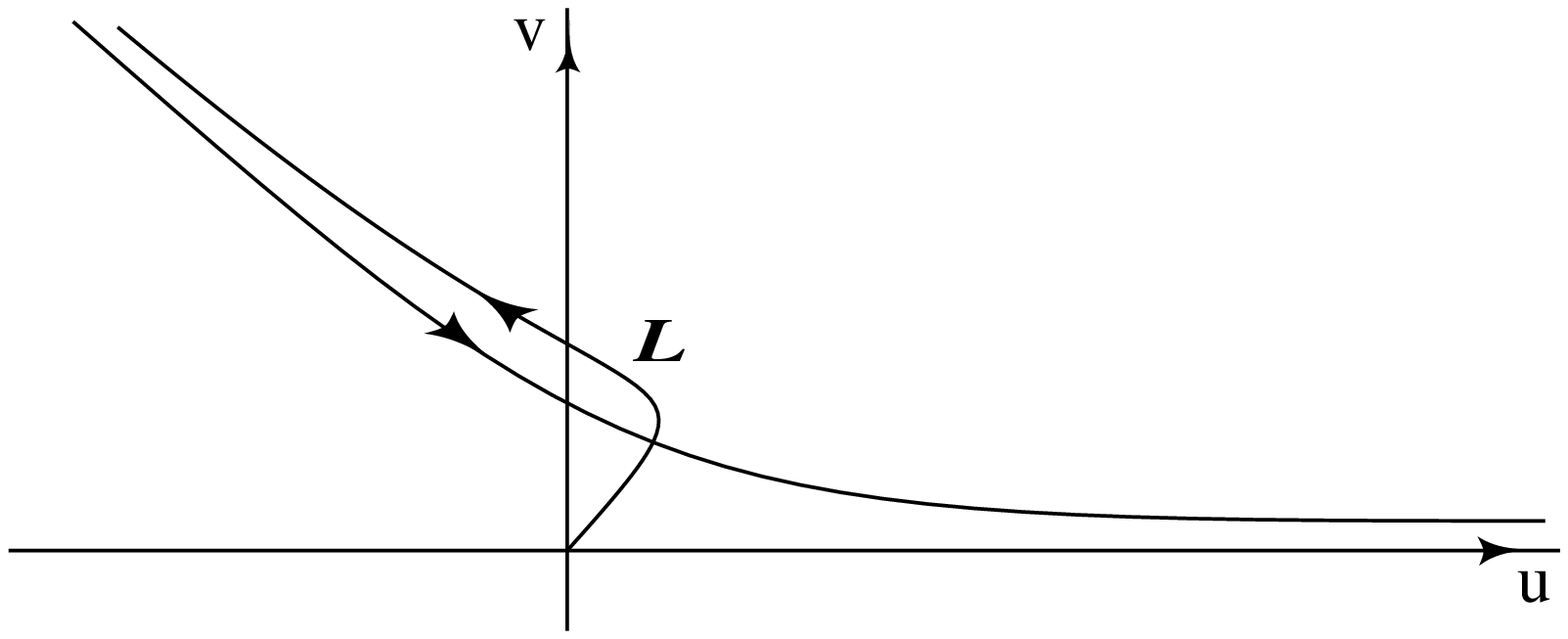,width=10.cm}}}	
{\bf Figure 3.2.}\quad The two contours for the integrals described in Remark
3.1. 
\vspace*{0.4cm}

\remark
Comparing the representations of $\Hi(z)$ in (\ref{d6}) and (\ref{d12}) 
we observe that the Airy function in (\ref{d12})
 disappears  as $z$ crosses the half line 
$\phase\ z=2\pi/3$. On that line the  argument of the Airy function 
$\Ai\left(ze^{-2\pi i/3}\right)$ becomes positive (see (\ref{d9})).
 Thus the dominance of $\Hi(z)$ over 
$\Ai\left(ze^{-2\pi i/3}\right)$ is maximal at this line,
this is therefore the place where the exponentially-small contribution is
``switched on''. 
This aspect is connected with the Stokes phenomenon
in asymptotics,  and the half-line $\phase\ z=2\pi/3$ is 
a Stokes line for $\Hi(z)$; see \cite{Par}. 

\subsection{The Scorer function $\Gi(z)$}%
It is convenient to have a direct method for $\Gi(z)$ that 
is not based on results for the $\Hi-$function, especially 
if $z$ is positive. We give only a few details on this case. The method 
can be used for the sector $0\le\phase\ z\le2\pi/3$. 

The first step is 
to replace the sine function by two exponentials. Then we obtain 
\begin{equation}
\Gi(z)=\frac1{2\pi i}\left[\Gi^+(z)-\Gi^-(z)\right],
\quad \Gi^{\pm}(z)=\int_{0}^{\infty} e^{\pm i\psi(t)}\,dt,\label{d13}
\end{equation}
where $\psi(t)=zt+\frac13t^3=\psi_r(u,v)+i\psi_i(u,v)$, with

\begin{equation}
\begin{array}{ll}
\psi_r(u,v)&=\frac13u^3-uv^2+xu-yv,\\
\psi_i(u,v)&=u^2v-\frac13v^3+xv+yu.\label{d14}
\end{array}
\end{equation}

The path of integration for $\Gi^\pm(z)$ is found by 
solving the equation $\psi_r(u,v)=0$. 

\noindent
For $\Gi^+(z)$ the path runs into the valley at
$\infty\exp(\pi i/6)$. The result is
\begin{equation}
\Gi^+(z)=\int_{0}^{\infty} e^{-\psi_i(u,v)}\,g(u)\,du, \quad
g(u)=1+i\frac{dv}{du}=1+i\frac{u^2-v^2+x}{2uv+y}.\label{d15}
\end{equation} 
For  $\Gi^-(z)$ the integral along $[0,\infty)$ can be replaced to 
a path along the half-line with $\phase\ t=-\pi/6$. In this valley 
no real solution of $\psi_r(u,v)=0$ is available, and we take a path 
that first runs into the valley at 
$\infty\exp(-5\pi i/6)$ and then returns to the valley at 
$\infty\exp(-\pi i/6)$. This introduces an Airy function, and we obtain
\begin{equation}
\Gi^-(z)=-\int_{0}^{\infty} e^{-\psi_i(u,v)}\,g(u)\,du+2\pi\Ai(z),\label{d16}
\end{equation}
where $g(u)$ is as in (\ref{d15}).
Adding the results in (\ref{d15}) and (\ref{d16}), we obtain a simple
non-oscillating integral plus an Airy function:
\begin{equation}
\Gi(z)=\frac1{\pi i}\int_{0}^{\infty} e^{-
\psi_i(u,v)}\,g(u)\,du+i\Ai(z),\label{d17}
\end{equation}
where the relation between $u$ and $v$ is given by
$$v=\frac{-y+\sqrt{{y^2+4u^2(x+\frac13u^2)}}}{2u}.$$

If $z$ is real and non-negative (\ref{d17}) becomes real.
 The term with the Airy
function is canceled by the imaginary contribution of $g(u)/i$. 
The remaining integral in
(\ref{d17}) should be modified in this case. The contour runs from
 the origin to the saddle point  $i\sqrt{{x}}$, and from this
 point into the valley at
$\infty\exp(\pi i/6)$. Integrating with respect to $v$, we obtain the
real representation
\begin{equation}
\Gi(x)=\frac1{\pi}\left[\int_0^{\sqrt{{x}}}e^{-
xv+\frac13v^3}\,dv+\int_{\sqrt{{x}}}^\infty
e^{2xv-\frac83v^3}\,dv\right],\quad x\ge0.\label{d18}
\end{equation}

\section{Numerical illustrations}

 We give some numerical results which serve
 as demonstration of our method. 
 
 In order to evaluate the Scorer functions in the whole complex plane, we
 need to select software for the evaluation of the Airy functions of complex
 arguments and for
 the quadrature of real functions 
 over an infinite integral. For the first purpose, we
 use the public domain  
 subroutines ZAIRY and ZBIRY by D.E. Amos \cite{Amo} 
 and for the semi-infinite integral we use the automatic adaptative 
 integrator DQAGIE by R. Piessens. All these codes can be retrieved from
 the SLATEC public domain library \cite{Sla}
  (see also GAMS: guide to available
 mathematical software \cite{Gam}).

 The connection formulae given in Section 2 together with the non-oscillating
 integrals derived in Section 3 can be used to evaluate $\Hi (z)$ and
 $\Gi (z)$ in the whole complex plane.  
 By using the integral 
 representations for $\Hi (z)$ in the domain $\pi \le \phase\
  z\le 2\pi/3$ the 
 following stable algorithm can be considered: 
 
 \begin{algorithm} {(via (\ref{d6}))}
 \begin{itemize}
 \item{} Whenever $\mbox{Im}\ (z)<0$, use eq. (\ref{t9}).
 \item{} If $z\in S_{1}^{(2)}$ then obtain $\Hi (z)$ by quadrature.
 \item{} If $z\in S_{1}^{(1)}\bigcup S_{0}$ obtain $\Hi (z)$ via (\ref{t7}).
 \item{} Obtain $\Gi (z)$ everywhere in the complex plane by using (\ref{t8}). 
 \end{itemize}
 \label{al1}
 \end{algorithm}

\noindent
 where $S_{1}^{(1)}$ is the sector $\pi /3 \le \mbox{ph}\, z\le 2\pi /3$ and 
$S_{1}^{(2)}$ is the sector $2\pi /3 \le \mbox{ph}\le  \pi$.

  However, in the fourth step two integrals for $\Hi$ will be needed. 
 Thus, the following 
 stable scheme is expected to be more efficient provided fast algorithms
 to compute the Airy 
 functions $\Ai$ and $\Bi$ are available:

\newpage
 \begin{algorithm} {(via (\ref{d6}) and (\ref{d17}))}
 \begin{itemize}
 \item{} Whenever $\mbox{Im}\ (z)<0$, use eq. (\ref{t9}).
 \item{} If $z\in S_{1}^{(2)}$ then obtain $\Hi (z)$ by quadrature.
 \item{} If $z\in S_{1}^{(1)}\bigcup S_{0}$ then obtain $\Gi (z)$ by 
 quadrature via  (\ref{d17}) or (\ref{d18}).
 \item{} For the remaining
 cases, apply eq. (\ref{i6}).
 \end{itemize}
 \label{al2}
 \end{algorithm}
 
  The second algorithm is preferred in most circumstances.
 However, we have experienced
 that the computation of $\Gi(z)$ when $\mbox{Re}(z) >0$ 
 and $\mbox{Im}(z)\rightarrow 0$ is more
 efficient when the first algorithm is considered. 
 Probably, the best numerical 
 strategy is a combination of both algorithms, 
 together with the use of asymptotic expansions
 for large $|z|$ and series expansions for small $|z|$. 
 The best strategy may also depend on the choice of the quadrature
 rule.
 This numerical study lies beyond the
 scope of the present paper.

  We end this section by showing numerical results. 
 Of interest are the sectors where (\ref{t5}) and (\ref{t6})
 are valid, namely, 
 $S_{1}^{(2)}$ for $\Hi (z)$ and
  $S_{1}^{(1)}\bigcup \mathcal{S}_{0}$ for $\Gi (z)$.
  Of particular interest are the regions 
  $S_{0}$ 
  for $\Gi (z)$ and $S_{1}^{(2)}$ for $\Hi (z)$
  since in this case we compute the functions
  directly by quadrature and the corresponding
  integrals can be compared with
  asymptotics if $z$ is large. We are giving explicit results with
  an accuracy of 8
  digits for $\Hi (z)$ which are compared with the asymptotic
  expantion (\ref{t5}) up to order $1/z^{10}$ 
  whenever this is possible. 
 
 The results in Table 4.1 are 
obtained by means of a Fortran program coded in double
 precision arithmetic in which
 the integral (\ref{d6}) is evaluated. Together with the
 results, we show the number of integration
 steps needed to attain an accuracy of 8 digits. 
 The results from asymptotics (shown inside brackets) are
 seen to coincide with those from the integral (\ref{d6}) for
 $|z|=100$, but for $|z|=10$ we observe discrepancies in the last digits
 which are due to the failure of the asymptotic expansion ($z$ is not
 large enough).
 We also found agreement with the asymptotic expansion
 for $\Hi (z)$ for large $z$ in the sector 
 $S_{1}^{(2)}$.
 For real negative $x$ our results coincide with those
 given by Scorer \cite{Sco}.

 One sees that at the Stokes line ph $z=2\pi /3$ 
the quadrature requires more steps, as can be expected
 given the appearance of a discontinuity in $dv/du$
 at the maximum $v$. On the other hand,
  the faster convergence takes place
 when we are far from the Stokes line.
  Moreover, convergence tends to be slower as $|z|$ becomes
  smaller; as $|z|$ becomes larger
 the effect of the singularity in the derivative $dv/du$ for 
 $\phase\ z=\frac23 \pi$ would 
 appear at larger $u$, where the exponential in the
 integrand is smaller.

  Similarly, one can test the performance of the integral
 representation (plus an Airy function) for $\Gi (z)$ 
 in the sector $ 0 \le$ ph $z \le 2\pi /3$, which should coincide
 with the results from the asymptotic expansion (\ref{t6})
 for $|\phase\ z|<\pi /3$ and large $z$.
 We also find agreement with the results tabulated by Scorer
 for real positive $x$.
 One observes that the convergence of
 the integral representation (\ref{d17})  becomes slower as
 we approach the  real line $\phase\ z=0$. 
 With the quadrature DQAGIE we obtain a better performance using
 Algorithm \ref{al1} in this case.

 \newpage

  \vspace*{0.3cm} 
  \begin{tabular}{|l|l|l|l|l|}
  \hline
      &     & $\phase\ z=\pi$  & $\phase\ z=5\pi/6$& $\phase\ z=2\pi /3$ \\
  \hline
  $|z|=1$ & Re ($\Hi$)& 0.22066961      & 0.22331566    & 0.23477589     \\
      &             &    $\{195\}$ &    $\{195\}$ &    $\{345\}$  \\     
  \hline
  $|z|=1$ & Im ($\Hi$)& 0       & 6.2133021 $10^{-2}$   & 0.13605894      \\
      &             &         & $\{165\}$   &  $\{465\}$      \\
  \hline 
  \hline
  $|z|=10$&Re ($\Hi$) & 3.1768535 $10^{-2}$ & 2.7597145 $10^{-2}$ &
   1.5948003 $10^{-2}$    \\
       &  &  $\{75\}$  &  $\{75\}$       &   $\{135\}$  \\     
    & & (3.1768528 $10^{-2}$) & (2.7597137 $10^{-2}$) &
   (1.5947998 $10^{-2}$)            \\ 
  \hline 
  $|z|=10$ & Im ($\Hi$)& 0 & 1.5859789 $10^{-2}$ & 2.7622751   $10^{-2}$     \\
       &  &  &  $\{75\}$ &     $\{225\}$  \\  
       & & &  (1.5859786 $10^{-2}$) & (2.7622742 $10^{-2}$)           \\
  \hline
  \hline
   $|z|=100$ & Re ($\Hi$)& 3.1830925 $10^{-3}$ & 2.7566477 $10^{-3}$ 
   & 1.5915526 $10^{-3}$  \\
   & & $\{135\}$  &  $\{165\}$ &   $\{165\}$ \\        
  & & (3.1830925 $10^{-3}$) & (2.7566477 $10^{-3}$) & (1.5915526 $10^{-3}$) \\ 
  \hline 
  $|z|=100$ & Im ($\Hi$)& 0  & 1.5915439 $10^{-3}$   & 2.7566500 $10^{-3}$  \\
        &             &    & $\{165\}$             & $\{165\}$             \\
       &              &    & (1.5915439 $10^{-3}$)& (2.7566500 $10^{-3}$)  \\
 \hline 
 \end{tabular}

 \vspace*{0.2cm}
 \noindent
 {\bf Table 4.1} The real and imaginary parts of $\Hi (z)$ in the sector
 $\pi \le \phase\ z \le \frac{2}{3} \pi$. The result from
 the asymptotic expansion up to order $1/z^{10}$ is shown inside brackets. 
 The number of integration steps
 for each evaluation is shown within the braces.

 \vspace*{0.3cm}
   
 As a further illustration, we will give two plots for $|z|=1$
 obtained by using Algorithm \ref{al2}.
 The graphs 
 show the real and imaginary parts of the Scorer functions.
 One observes the smooth connection
 between the different sectors in the complex plane. 

\vspace*{0.3cm}

\begin{minipage}{6.cm}
\centerline{\protect\hbox{\psfig{file=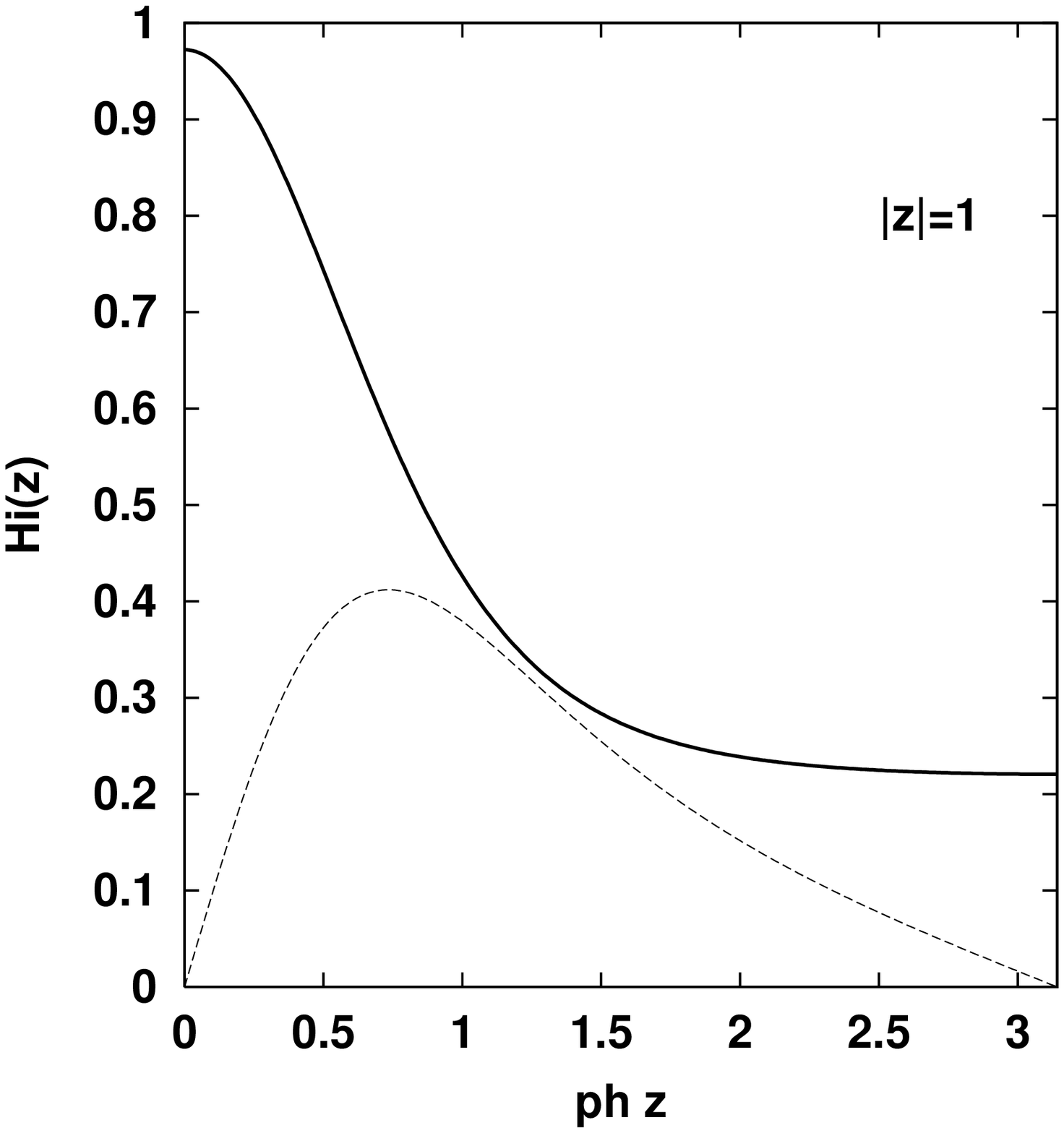,width=6.2cm}}}
\end{minipage}
\hfill
\begin{minipage}{6.cm} 
\centerline{\protect\hbox{\psfig{file=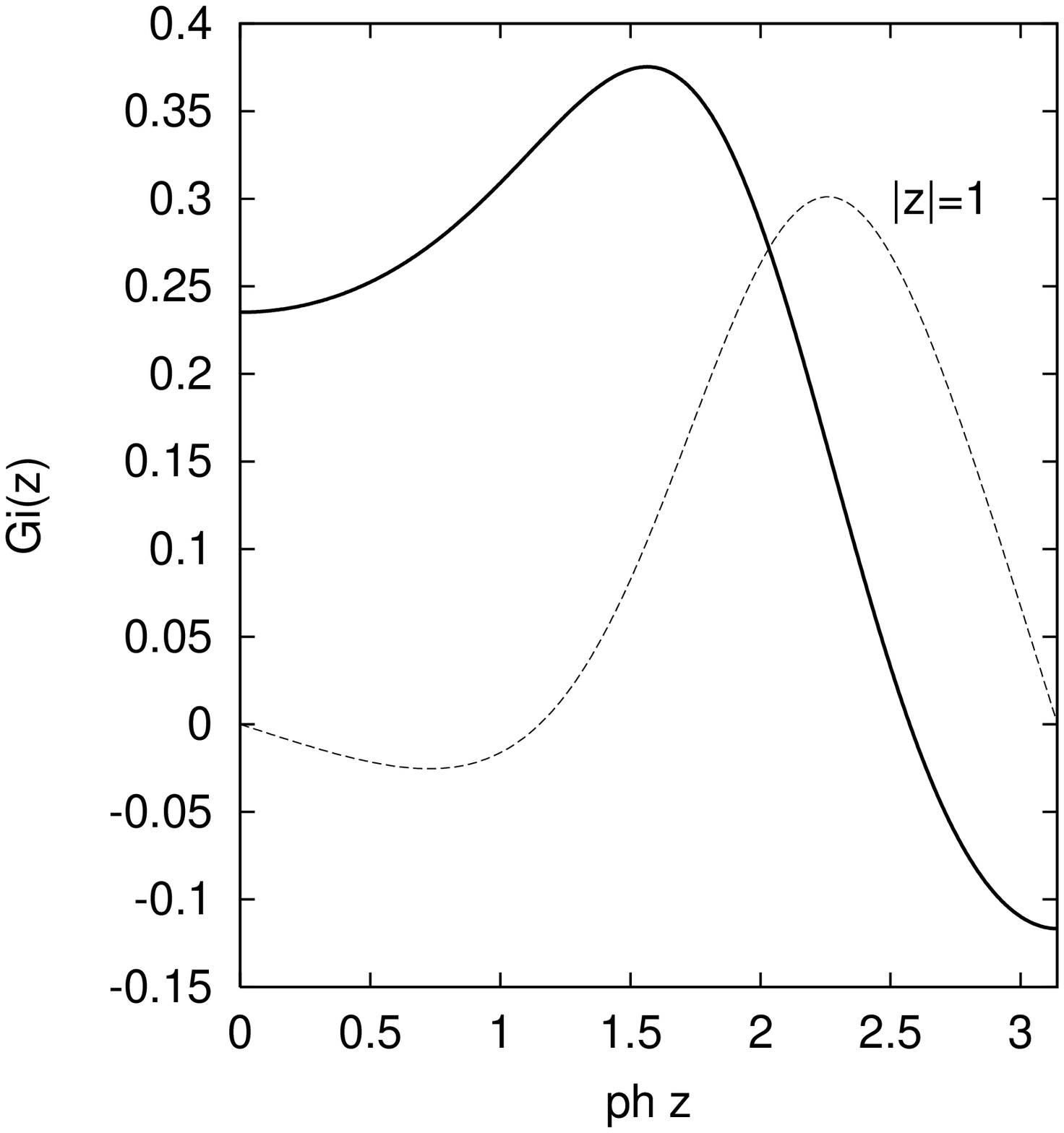,width=6.2cm}}}
\end{minipage}

\vspace*{0.2cm}
\noindent
{\bf Figure 4.2} The real (solid) and imaginary (dashed) 
parts of $\Hi (z)$ (left) and $\Gi (z)$ (right) 
for $|z|=1$ and $0\le \phase\ z \le \pi$. 

\vspace*{0.4cm}

\section{Summary and concluding remarks}
As mentioned in the cited references (see \cite{Lee} and \cite{Leo})
the inhomogeneous Airy functions (Scorer functions) are used in 
several physical problems. We have used functional
relations and derived integral representations of the Scorer functions 
 that can be used for stable numerical computations for
all complex values of the argument $z$. The integrals can be easily split up
into real and imaginary parts.

For the Scorer functions Maclaurin and asymptotic expansions are
available. To bridge the gap between the areas where convergent or 
asymptotic series can be used, one can use the representations in this
paper, although they can be used for all values of the argument. The tool
one needs is a suitable quadrature method for computing real integrals on
unbounded real intervals that converge very fast at infinity. We have
illustrated the method by giving a few numerical results based on
selecting a quadrature rule; we have not investigated an optimal
choice of quadrature rule for computing the Scorer functions.

We have shown how to handle oscillating integrals for a certain 
set of special functions, and this is quite instructive for applying the
method  to other functions. A similar method can be used for 
the Airy functions and another treatment can be found in \cite{Tem1}
for modified Bessel functions of imaginary order.
Still many special functions need to be considered in order to get 
reliable software, in particular for complex and/or large values of 
the parameters.

\bigskip
\noindent
{\bf Acknowledgments.\  }
 A.G. and J.S. would like to acknowledge the hospitality of CWI during
     their stay. A.G. and J.S. also acknowledge financial support from
     the Conseller\'{\i}a de Educaci\'on y Ciencia (Generalitat Valenciana).

\noindent
The authors thank the referee for valuable suggestions.

\bibliographystyle{amsplain}

\begin{thebibliography}{99}

\bibitem{Abr} M. Abramowitz and I.A. Stegun (Eds.),
\textit{Handbook of Mathematical functions},
National Bureau of Standards Applied Mathematics Series No. 55. U.S.
Government Printing Office, Washington, DC.  

\bibitem{Amo} D.E. Amos. 
``Algorithm 644: A portable package for Bessel functions of 
a complex argument and nonnegative order''.
 ACM Trans. Math. Softw. 12 (1986) 265-273.

\bibitem{Cor} R.M. Corless, D.J. Jeffrey and H. Rasmussen
``Numerical evaluation of Airy functions with complex arguments''.
J. Comput. Phys. 99 (1992), 106-114.

\bibitem{Ext} H. Exton.
``The asymptotic behaviour of the inhomogeneous Airy function $\Hi(z)$''. 
Math. Chronicle 12 (1983),99-104.

\bibitem{Fab} B. Fabijonas
``The computation of Scorer functions''. 
Lecture during the 1998 Annual SIAM Meeting in Toronto, Canada. 

\bibitem{Gam} GAMS: Guide to available mathematical software.
http://gams.nist.gov

\bibitem{Lee} S.-Y. Lee
``The inhomogeneous Airy functions, $\Gi(z)$ and $\Hi(z)$''.
 J. Chem. Phys. 72 (1980), 332-336.

\bibitem{Loz} D.W. Lozier and F.W.J. Olver.
 ``Numerical evaluation of special functions''. In W. Gautschi (Ed.),
 AMS Proceedings of Symposia in Applied Mathematics
 48 (1998), pp. 79--125.

\bibitem{Leo} A.J. MacLeod.
``Computation of inhomogeneous Airy functions''.
J. Comput. Appl. Math. 53 (1994) 109-116.

\bibitem{Nis} The National Institute of Standards and Technology
has a public web site that includes an extensive 
treatment of Scorer functions: http://www.nist.gov/DigitalMathLib.

\bibitem{Olv} F.W.J. Olver.
\textit{Asymptotics and Special Functions}.
Academic Press, New York. Reprinted in 1997 by A.K. Peters.

\bibitem{Par} R.B. Paris and A.D. Wood. 
 ``Stokes phenomenon demystified'',
 IMA Bulletin 31 (1995) No.1-2,21-28.

\bibitem{Sla} 
 SLATEC Public Domain Mathematical Library. 

 gopher://archives.math.utk.edu/11/software/multi-platform/SLATEC


\bibitem{Sch} Z. Schulten, D.G.M. Anderson, and R.G. Gordon.
``An algorithm for the evaluation of complex Airy functions''.
 J. Comput. Phys. 31 (1979) 60-75. 

\bibitem{Sco} R.S. Scorer. 
``Numerical evaluation of integrals of the form 
$I=\int_{x_1}^{x_2}\,f(x)e^{i\phi(x)}\,dx$
and the tabulation of the function 
$\Gi(z)=(1/\pi) \int_{0}^{\infty}\,\sin\left(uz+\frac{1}{3}u^3\right)\,du$''.
Quart. J. Mech. Appl. Math. 3 (1950) 107-112. 

\bibitem{Tem1} N.M. Temme.
``Steepest descent paths for integrals defining the modified 
Bessel functions of imaginary order''. 
Methods Appl. Anal. 1 (1994) 14--24.

\bibitem{Tem2} N.M. Temme. 
\textit{Special functions: An introduction to the classical functions
of mathematical physics}.
John Wiley \& Sons, New York, 1996. 

\bibitem{Won} R. Wong
\textit{Asymptotic approximations of integrals}.
Academic Press, New York, 1989.

\end{thebibliography}

\end{document}